
\magnification 1200
\input amssym.def
\input amssym
\catcode`\@=11
\font\eightrm=cmr8         \font\eighti=cmmi8
\font\eightsy=cmsy8        \font\eightbf=cmbx8
\font\eighttt=cmtt8        \font\eightit=cmti8
\font\eightsl=cmsl8        \font\sixrm=cmr6
\font\sixi=cmmi6           \font\sixsy=cmsy6
\font\sixbf=cmbx6
\font\tengoth=eufm10       \font\tenbboard=msbm10
\font\eightgoth=eufm8      \font\eightbboard=msbm8
\font\sevengoth=eufm7      \font\sevenbboard=msbm7
\font\sixgoth=eufm6        \font\fivegoth=eufm5
\font\tencyr=wncyr10       
\font\eightcyr=wncyr8      
\font\sevencyr=wncyr7      
\font\sixcyr=wncyr6
\newfam\gothfam           \newfam\bboardfam
\newfam\cyrfam
\def\tenpoint{%
  \textfont0=\tenrm \scriptfont0=\sevenrm \scriptscriptfont0=\fiverm
  \def\rm{\fam\z@\tenrm}%
  \textfont1=\teni  \scriptfont1=\seveni  \scriptscriptfont1=\fivei
  \def\oldstyle{\fam\@ne\teni}\let\old=\oldstyle
  \textfont2=\tensy \scriptfont2=\sevensy \scriptscriptfont2=\fivesy
  \textfont\gothfam=\tengoth \scriptfont\gothfam=\sevengoth
  \scriptscriptfont\gothfam=\fivegoth
  \def\goth{\fam\gothfam\tengoth}%
  \textfont\bboardfam=\tenbboard \scriptfont\bboardfam=\sevenbboard
  \scriptscriptfont\bboardfam=\sevenbboard
  \def\bb{\fam\bboardfam\tenbboard}%
 \textfont\cyrfam=\tencyr \scriptfont\cyrfam=\sevencyr
  \scriptscriptfont\cyrfam=\sixcyr
  \def\cyr{\fam\cyrfam\tencyr}%
  \textfont\itfam=\tenit
  \def\it{\fam\itfam\tenit}%
  \textfont\slfam=\tensl
  \def\sl{\fam\slfam\tensl}%
  \textfont\bffam=\tenbf \scriptfont\bffam=\sevenbf
  \scriptscriptfont\bffam=\fivebf
  \def\bf{\fam\bffam\tenbf}%
  \textfont\ttfam=\tentt
  \def\tt{\fam\ttfam\tentt}%
  \abovedisplayskip=9pt plus 3pt minus 9pt
  \belowdisplayskip=\abovedisplayskip
  \abovedisplayshortskip=0pt plus 3pt
  \belowdisplayshortskip=4pt plus 3pt
  \smallskipamount=3pt plus 1pt minus 1pt
  \medskipamount=6pt plus 2pt minus 2pt
  \bigskipamount=12pt plus 4pt minus 4pt
  \normalbaselineskip=12pt
  \setbox\strutbox=\hbox{\vrule height8.5pt depth3.5pt width0pt}%
  \let\bigf@nt=\tenrm       \let\smallf@nt=\sevenrm
  \normalbaselines\rm}

\def\eightpoint{%
  \textfont0=\eightrm \scriptfont0=\sixrm \scriptscriptfont0=\fiverm
  \def\rm{\fam\z@\eightrm}%
  \textfont1=\eighti  \scriptfont1=\sixi  \scriptscriptfont1=\fivei
  \def\oldstyle{\fam\@ne\eighti}\let\old=\oldstyle
  \textfont2=\eightsy \scriptfont2=\sixsy \scriptscriptfont2=\fivesy
  \textfont\gothfam=\eightgoth \scriptfont\gothfam=\sixgoth
  \scriptscriptfont\gothfam=\fivegoth
  \def\goth{\fam\gothfam\eightgoth}%
  \textfont\cyrfam=\eightcyr \scriptfont\cyrfam=\sixcyr
  \scriptscriptfont\cyrfam=\sixcyr
  \def\cyr{\fam\cyrfam\eightcyr}%
  \textfont\bboardfam=\eightbboard \scriptfont\bboardfam=\sevenbboard
  \scriptscriptfont\bboardfam=\sevenbboard
  \def\bb{\fam\bboardfam}%
  \textfont\itfam=\eightit
  \def\it{\fam\itfam\eightit}%
  \textfont\slfam=\eightsl
  \def\sl{\fam\slfam\eightsl}%
  \textfont\bffam=\eightbf \scriptfont\bffam=\sixbf
  \scriptscriptfont\bffam=\fivebf
  \def\bf{\fam\bffam\eightbf}%
  \textfont\ttfam=\eighttt
  \def\tt{\fam\ttfam\eighttt}%
  \abovedisplayskip=9pt plus 3pt minus 9pt
  \belowdisplayskip=\abovedisplayskip
  \abovedisplayshortskip=0pt plus 3pt
  \belowdisplayshortskip=3pt plus 3pt
  \smallskipamount=2pt plus 1pt minus 1pt
  \medskipamount=4pt plus 2pt minus 1pt
  \bigskipamount=9pt plus 3pt minus 3pt
  \normalbaselineskip=9pt
  \setbox\strutbox=\hbox{\vrule height7pt depth2pt width0pt}%
  \let\bigf@nt=\eightrm     \let\smallf@nt=\sixrm
  \normalbaselines\rm}

\tenpoint
\def\pc#1{\bigf@nt#1\smallf@nt}         \def\pd#1 {{\pc#1} }

\def\Grille{\setbox13=\vbox to 5\unitlength{\hrule width 109mm\vfill}
\setbox13=\vbox to 65mm{\offinterlineskip\leaders\copy13\vfill\kern-1pt\hrule}
\setbox14=\hbox to 5\unitlength{\vrule height 65mm\hfill}
\setbox14=\hbox to 109mm{\leaders\copy14\hfill\kern-2mm\vrule height 65mm}
\ht14=0pt\dp14=0pt\wd14=0pt \setbox13=\vbox to
0pt{\vss\box13\offinterlineskip\box14} \wd13=0pt\box13}


\def\fleche(#1,#2)\dir(#3,#4)\long#5{%
\noalign{\leftput(#1,#2){\vector(#3,#4){#5}}}}

\def\ligne(#1,#2)\dir(#3,#4)\long#5{%
\noalign{\leftput(#1,#2){\lline(#3,#4){#5}}}}

\def\put(#1,#2)#3{\noalign{\setbox1=\hbox{%
    \kern #1\unitlength
    \raise #2\unitlength\hbox{$#3$}}%
    \ht1=0pt \wd1=0pt \dp1=0pt\box1}}


\def\diagram#1{\def\normalbaselines{\baselineskip=0pt\lineskip=5pt}
\matrix{#1}}

\def\gfl#1#2#3{\smash{\mathop{\hbox to#3{\leftarrowfill}}\limits
^{\scriptstyle#1}_{\scriptstyle#2}}}

\def\vfl#1#2#3{\llap{$\scriptstyle #1$}
\left\downarrow\vbox to#3{}\right.\rlap{$\scriptstyle #2$}}


 \message{`lline' & `vector' macros from LaTeX}
 \catcode`@=11
\def\{{\relax\ifmmode\lbrace\else$\lbrace$\fi}
\def\}{\relax\ifmmode\rbrace\else$\rbrace$\fi}
\def\newcount{\alloc@0\count\countdef\insc@unt}
\def\newdimen{\alloc@1\dimen\dimendef\insc@unt}
\def\newwrite{\alloc@7\write\chardef\sixt@@n}

\newwrite\@unused
\newcount\@tempcnta
\newcount\@tempcntb
\newdimen\@tempdima
\newdimen\@tempdimb
\newbox\@tempboxa

\def\@spaces{\space\space\space\space}
\def\@whilenoop#1{}
\def\@whiledim#1\do #2{\ifdim #1\relax#2\@iwhiledim{#1\relax#2}\fi}
\def\@iwhiledim#1{\ifdim #1\let\@nextwhile=\@iwhiledim
        \else\let\@nextwhile=\@whilenoop\fi\@nextwhile{#1}}
\def\@badlinearg{\@latexerr{Bad \string\line\space or \string\vector
   \space argument}}
\def\@latexerr#1#2{\begingroup
\edef\@tempc{#2}\expandafter\errhelp\expandafter{\@tempc}%
\def\@eha{Your command was ignored.
^^JType \space I <command> <return> \space to replace it
  with another command,^^Jor \space <return> \space to continue without
it.}
\def\@ehb{You've lost some text. \space \@ehc}
\def\@ehc{Try typing \space <return>
  \space to proceed.^^JIf that doesn't work, type \space X <return> \space to
  quit.}
\def\@ehd{You're in trouble here.  \space\@ehc}

\typeout{LaTeX error. \space See LaTeX manual for explanation.^^J
 \space\@spaces\@spaces\@spaces Type \space H <return> \space for
 immediate help.}\errmessage{#1}\endgroup}
\def\typeout#1{{\let\protect\string\immediate\write\@unused{#1}}}

\font\tenln    = line10
\font\tenlnw   = linew10
\newdimen\@wholewidth
\newdimen\@halfwidth
\newdimen\unitlength

\unitlength =1pt


\def\thinlines{\let\@linefnt\tenln \let\@circlefnt\tencirc
  \@wholewidth\fontdimen8\tenln \@halfwidth .5\@wholewidth}
\def\thicklines{\let\@linefnt\tenlnw \let\@circlefnt\tencircw
  \@wholewidth\fontdimen8\tenlnw \@halfwidth .5\@wholewidth}

\def\linethickness#1{\@wholewidth #1\relax \@halfwidth .5\@wholewidth}

\newif\if@negarg

\def\lline(#1,#2)#3{\@xarg #1\relax \@yarg #2\relax
\@linelen=#3\unitlength
\ifnum\@xarg =0 \@vline
  \else \ifnum\@yarg =0 \@hline \else \@sline\fi
\fi}

\def\@sline{\ifnum\@xarg< 0 \@negargtrue \@xarg -\@xarg \@yyarg -\@yarg
  \else \@negargfalse \@yyarg \@yarg \fi
\ifnum \@yyarg >0 \@tempcnta\@yyarg \else \@tempcnta -\@yyarg \fi
\ifnum\@tempcnta>6 \@badlinearg\@tempcnta0 \fi
\setbox\@linechar\hbox{\@linefnt\@getlinechar(\@xarg,\@yyarg)}%
\ifnum \@yarg >0 \let\@upordown\raise \@clnht\z@
   \else\let\@upordown\lower \@clnht \ht\@linechar\fi
\@clnwd=\wd\@linechar
\if@negarg \hskip -\wd\@linechar \def\@tempa{\hskip -2\wd\@linechar}\else
     \let\@tempa\relax \fi
\@whiledim \@clnwd <\@linelen \do
  {\@upordown\@clnht\copy\@linechar
   \@tempa
   \advance\@clnht \ht\@linechar
   \advance\@clnwd \wd\@linechar}%
\advance\@clnht -\ht\@linechar
\advance\@clnwd -\wd\@linechar
\@tempdima\@linelen\advance\@tempdima -\@clnwd
\@tempdimb\@tempdima\advance\@tempdimb -\wd\@linechar
\if@negarg \hskip -\@tempdimb \else \hskip \@tempdimb \fi
\multiply\@tempdima \@m
\@tempcnta \@tempdima \@tempdima \wd\@linechar \divide\@tempcnta \@tempdima
\@tempdima \ht\@linechar \multiply\@tempdima \@tempcnta
\divide\@tempdima \@m
\advance\@clnht \@tempdima
\ifdim \@linelen <\wd\@linechar
   \hskip \wd\@linechar
  \else\@upordown\@clnht\copy\@linechar\fi}

\def\@hline{\ifnum \@xarg <0 \hskip -\@linelen \fi
\vrule height \@halfwidth depth \@halfwidth width \@linelen
\ifnum \@xarg <0 \hskip -\@linelen \fi}

\def\@getlinechar(#1,#2){\@tempcnta#1\relax\multiply\@tempcnta 8
\advance\@tempcnta -9 \ifnum #2>0 \advance\@tempcnta #2\relax\else
\advance\@tempcnta -#2\relax\advance\@tempcnta 64 \fi
\char\@tempcnta}

\def\vector(#1,#2)#3{\@xarg #1\relax \@yarg #2\relax
\@linelen=#3\unitlength
\ifnum\@xarg =0 \@vvector
  \else \ifnum\@yarg =0 \@hvector \else \@svector\fi
\fi}

\def\@hvector{\@hline\hbox to 0pt{\@linefnt
\ifnum \@xarg <0 \@getlarrow(1,0)\hss\else
    \hss\@getrarrow(1,0)\fi}}

\def\@vvector{\ifnum \@yarg <0 \@downvector \else \@upvector \fi}

\def\@svector{\@sline
\@tempcnta\@yarg \ifnum\@tempcnta <0 \@tempcnta=-\@tempcnta\fi
\ifnum\@tempcnta <5
  \hskip -\wd\@linechar
  \@upordown\@clnht \hbox{\@linefnt  \if@negarg
  \@getlarrow(\@xarg,\@yyarg) \else \@getrarrow(\@xarg,\@yyarg) \fi}%
\else\@badlinearg\fi}

\def\@getlarrow(#1,#2){\ifnum #2 =\z@ \@tempcnta='33\else
\@tempcnta=#1\relax\multiply\@tempcnta \sixt@@n \advance\@tempcnta
-9 \@tempcntb=#2\relax\multiply\@tempcntb \tw@
\ifnum \@tempcntb >0 \advance\@tempcnta \@tempcntb\relax
\else\advance\@tempcnta -\@tempcntb\advance\@tempcnta 64
\fi\fi\char\@tempcnta}

\def\@getrarrow(#1,#2){\@tempcntb=#2\relax
\ifnum\@tempcntb < 0 \@tempcntb=-\@tempcntb\relax\fi
\ifcase \@tempcntb\relax \@tempcnta='55 \or
\ifnum #1<3 \@tempcnta=#1\relax\multiply\@tempcnta
24 \advance\@tempcnta -6 \else \ifnum #1=3 \@tempcnta=49
\else\@tempcnta=58 \fi\fi\or
\ifnum #1<3 \@tempcnta=#1\relax\multiply\@tempcnta
24 \advance\@tempcnta -3 \else \@tempcnta=51\fi\or
\@tempcnta=#1\relax\multiply\@tempcnta
\sixt@@n \advance\@tempcnta -\tw@ \else
\@tempcnta=#1\relax\multiply\@tempcnta
\sixt@@n \advance\@tempcnta 7 \fi\ifnum #2<0 \advance\@tempcnta 64 \fi
\char\@tempcnta}

\def\@vline{\ifnum \@yarg <0 \@downline \else \@upline\fi}

\def\@upline{\hbox to \z@{\hskip -\@halfwidth \vrule
  width \@wholewidth height \@linelen depth \z@\hss}}

\def\@downline{\hbox to \z@{\hskip -\@halfwidth \vrule
  width \@wholewidth height \z@ depth \@linelen \hss}}

\def\@upvector{\@upline\setbox\@tempboxa\hbox{\@linefnt\char'66}\raise
     \@linelen \hbox to\z@{\lower \ht\@tempboxa\box\@tempboxa\hss}}

\def\@downvector{\@downline\lower \@linelen
      \hbox to \z@{\@linefnt\char'77\hss}}

\thinlines

\newcount\@xarg
\newcount\@yarg
\newcount\@yyarg
\newcount\@multicnt
\newdimen\@xdim
\newdimen\@ydim
\newbox\@linechar
\newdimen\@linelen
\newdimen\@clnwd
\newdimen\@clnht
\newdimen\@dashdim
\newbox\@dashbox
\newcount\@dashcnt


\newbox\tbox
\newbox\tboxa
\def\leftzer#1{\setbox\tbox=\hbox to 0pt{#1\hss}%
     \ht\tbox=0pt \dp\tbox=0pt \box\tbox}

\def\rightzer#1{\setbox\tbox=\hbox to 0pt{\hss #1}%
     \ht\tbox=0pt \dp\tbox=0pt \box\tbox}

\def\centerzer#1{\setbox\tbox=\hbox to 0pt{\hss #1\hss}%
     \ht\tbox=0pt \dp\tbox=0pt \box\tbox}
\def\image(#1,#2)#3{\vbox to #1{\offinterlineskip
    \vss #3 \vskip #2}}

\def\leftput(#1,#2)#3{\setbox\tboxa=\hbox{%
    \kern #1\unitlength
    \raise #2\unitlength\hbox{\leftzer{#3}}}%
    \ht\tboxa=0pt \wd\tboxa=0pt \dp\tboxa=0pt\box\tboxa}

\def\rightput(#1,#2)#3{\setbox\tboxa=\hbox{%
    \kern #1\unitlength
    \raise #2\unitlength\hbox{\rightzer{#3}}}%
    \ht\tboxa=0pt \wd\tboxa=0pt \dp\tboxa=0pt\box\tboxa}

\def\centerput(#1,#2)#3{\setbox\tboxa=\hbox{%
    \kern #1\unitlength
    \raise #2\unitlength\hbox{\centerzer{#3}}}%
    \ht\tboxa=0pt \wd\tboxa=0pt \dp\tboxa=0pt\box\tboxa}

\unitlength=1mm
\def\dash{\dasharrow}
\mathcode`A="7041 \mathcode`B="7042 \mathcode`C="7043 \mathcode`D="7044
\mathcode`E="7045 \mathcode`F="7046 \mathcode`G="7047 \mathcode`H="7048
\mathcode`I="7049 \mathcode`J="704A \mathcode`K="704B \mathcode`L="704C
\mathcode`M="704D \mathcode`N="704E \mathcode`O="704F \mathcode`P="7050
\mathcode`Q="7051 \mathcode`R="7052 \mathcode`S="7053 \mathcode`T="7054
\mathcode`U="7055 \mathcode`V="7056 \mathcode`W="7057 \mathcode`X="7058
\mathcode`Y="7059 \mathcode`Z="705A

\def\spacedmath#1{\def\packedmath##1${\bgroup\mathsurround=0pt ##1\egroup$}%
\mathsurround#1 \everymath={\packedmath}\everydisplay={\mathsurround=0pt }}

\def\nospacedmath{\mathsurround=0pt \everymath={}\everydisplay={} }
\parindent0cm
\spacedmath{2pt}
\vsize = 25.1truecm
\hsize = 16truecm
\hoffset = -.15truecm
\voffset = -.5truecm
\baselineskip15pt
\def\ind{\par\hskip 1cm\relax}
\def\moins{\mathrel{\hbox{\vrule height 3pt depth -2pt
width 6pt}}}
\def\rond{\kern 1pt{\scriptstyle\circ}\kern 1pt}
\def\epi{\rightarrow \kern -3mm\rightarrow }
\def\longepi{\longrightarrow \kern -5mm\longrightarrow }
\def\mono{\lhook\joinrel\mathrel{\longrightarrow}}

\def\Sp{\mathop{\rm Spec\,}\nolimits}
\def\iso{\mathrel{\mathop{\kern 0pt\longrightarrow }\limits^{\sim}}}
\def\vb{vector bundle }
\def\al{$k$\kern -1.5pt -algebr}

\def\Ker{\mathop{\rm Ker}\nolimits}
\def\Coker{\mathop{\rm Coker}}

\def\Pic{\mathop{\rm Pic}\nolimits}

\def\dim{\mathop{\rm dim}\nolimits}

\def\up#1{\raise 1ex\hbox{\smallf@nt#1}}
\def\note#1#2{\footnote{\parindent
0.4cm$^#1$}{\vtop{\eightpoint\baselineskip12pt\hsize15.5truecm
\noindent #2}}\parindent 0cm}
\overfullrule=0pt
\def\tx{\kern-1.5pt -}
\catcode`\@=12
\null\vskip0.1cm
\def\iso{\mathrel{\mathop{\kern 0pt\longrightarrow }\limits^{\sim}}}
\def\vb{vector bundl}
\def\lb{line bundl}
\centerline{\bf Vector bundles on curves and generalized theta functions:}
\centerline{\bf recent results and open problems}
\smallskip \centerline{Arnaud {\pc BEAUVILLE} \footnote{\parindent
0.6cm (*)}{\vtop{\eightpoint\baselineskip12pt\hsize15.2truecm\noindent
 Partially supported by the
European Science project ``Geometry of Algebraic Varieties", Contract
 SCI-0398-C(A).}}} \vskip0.9cm

{\bf Introduction}
\ind  It is known essentially since
Riemann that one can associate to any compact Riemann surface $X$ an
Abelian variety, the Jacobian $JX$, together with a divisor $\Theta$
(well-defined up to translation) which can be defined both in a geometric way
and as the zero locus of an explicit function, the Riemann theta function. The
geometry of the pair $(JX,\Theta)$ is intricately (and beautifully) related to
the geometry of $X$.  \ind The idea that  higher rank vector bundles should
provide a non-Abelian analogue of the Jacobian appears already in the
influential paper [We] of A. Weil (though the notion of \vb e does not appear
as such in that paper!). The construction of the moduli spaces has been
achieved in the 60's, mainly by D.~Mumford  and the mathematicians of the
Tata Institute. However  it is only recently  that  the study of the
determinant line bundles on these moduli spaces  and of their spaces of
sections has made clear the analogy with the Jacobian. This is largely due to
the  intrusion  of Conformal Field Theory, where these spaces have appeared
(quite surprisingly for us!) as  fundamental objects.

\ind In these notes (based on a few lectures given in the Fall of {\old 1992}
at  MSRI, UCLA and University of Utah), I will try to give an overview of these
 new ideas.  I~must warn the reader that this is by no means intended to be a
complete account. I~have mainly focused on the determinant line bundles and
their
spaces of sections,  ignoring deliberately
 important areas like cohomology of the moduli spaces, moduli of Higgs
bundles, relations with integrable systems,  Langlands' geometric
correspondence..., simply because I felt it would have taken me too far afield.
For
the same reason I~haven't even tried to explain why the mathematical physicists
are so  interested in these moduli spaces.
 \vskip1cm  {\bf 1. The moduli space ${\cal SU}_X(r)$}
\smallskip
\ind Let $X$ be a compact Riemann surface of genus $g$. Recall that the
Jacobian
 $JX$ parametrizes line bundles of degree $0$ on $X$. We will also consider the
variety $J^{g-1}(X)$ which parametrizes line bundles of degree $g-1$ on $X$; it
carries a \ {\it canonical} Theta divisor\vskip-12pt
$$\Theta=\{M\in J^{g-1}(X)\ |\ H^0(X,M)\not=0\}\ .$$\vskip-3pt
For each line bundle $L$ on $X$ of degree $g-1$, the map $M\mapsto M\otimes
L^{-1}$ induces an isomorphism of
 $(J^{g-1}(X),\Theta)$ onto $(JX,\Theta_L)$, where  $\Theta_L$ is the divisor
on
$JX$ defined by $$\Theta_L=\{E\in JX \ |\ H^0(X,E\otimes L)\not=0\}\ .$$
\ind We know a great deal about the spaces $H^0(JX,{\cal O}(k\Theta))$. One of
the
key points is that the  sections of ${\cal O}(k\Theta)$ can be identified
with certain quasi-periodic functions on the universal cover of $JX$, the theta
functions of order $k$. In this way one gets for instance that the dimension of
$H^0(JX,{\cal O}(k\Theta))$ is $k^g$, that the
linear
system $|k\Theta|$ is  base-point free for $k\ge 2$ and  very ample  for
$k\ge 3$, and so on. \def\qplus{\mathop\oplus}
One even obtains a rather precise description of the
ring $\displaystyle  \qplus_{k\ge 0}H^0(JX,{\cal O}(k\Theta))$, the graded ring
of theta functions.\smallskip
 \ind The character who will play the role  of the
Jacobian in these lectures is the moduli  space ${\cal SU}_X(r)$ of
(semi-stable)
rank $r$ \vb es on $X$ with trivial determinant. It is an irreducible
projective
variety, whose points are isomorphism classes of \vb es which are direct sums
of
stable \vb es of degree $0$ (a degree $0$ \vb e $E$ is said to be {\it stable}
if
every proper subbundle of $E$ has degree $<0$).  By the theorem of Narasimhan
and
Seshadri, the points of ${\cal SU}_X(r)$ are also the  isomorphism classes of
representations $\pi_1(X)\longrightarrow SU(r)$ (hence the  notation ${\cal
SU}_X(r)$). The stable bundles form a smooth open subset of ${\cal SU}_X(r)$,
whose complement (which parametrizes decomposable bundles) is
singular\note{1}{Except in the cases $g\le 1$ and $g=r=2$, where the moduli
space
is smooth.}.  \ind The reason for fixing the determinant is that the moduli
space
${\cal U}_X(r)$ of \vb es of rank $r$ and degree $0$ is, up to a finite \'etale
covering, the product of ${\cal SU}_X(r)$ with $JX$, so the study of ${\cal
U}_X(r)$
is essentially reduced to that of ${\cal SU}_X(r)$. Of course the moduli spaces
${\cal SU}_X(r,L)$ of semi-stable bundles   \vb es
with a fixed determinant $L\in \Pic(X)$ is also of interest; for simplicity, in
these
lectures I~will concentrate on the most central  case $L={\cal O}_X$.  \ind
Observe
that when $g\le 1$ the spaces ${\cal SU}_X(r)$ consist only of direct sums of
line
bundles. Since these cases are quite easy to deal with directly,  I will
usually
assume implicitely  $g\ge 2$ in what follows.

 \vskip0.9cm{\bf 2. The determinant bundle}\vskip3pt
 \ind  The geometric definition of the theta  divisor extends in a natural way
to the
higher rank case. For any line bundle  $L\in J^{g-1}(X)$, define
$$\Theta_L=\{E\in {\cal SU}_X(r) \ |\ h^0(X,E\otimes L)\ge 1\}\ .$$
This turns out to be a Cartier divisor on ${\cal SU}_X(r)$ [D-N] (the key point
here
is that the degrees are chosen so that  $\chi(E\otimes L)=0$). The associated
line
bundle ${\cal L}:={\cal O}(\Theta_L)$  does not depend on the choice of $L$. It
is
called the {determinant bundle}, and will play a central role in our story.  It
is in
fact canonical, because of the following result (proved in [B 1] for
$r=2$ and in [D-N] in general):
\medskip
  {\pc THEOREM} 1. -- $\Pic {\cal
SU}_X(r)={\bf Z}{\cal L}$. \medskip
\ind By analogy with the rank one case, the global sections  of the line
bundles
${\cal L}^k$ are sometimes called {\it generalized theta functions} -- we will
briefly discuss this terminology in \S 10.
 \ind Can we describe $H^0( {\cal SU}_X(r), {\cal
L})$ and the map $\varphi^{}_{\cal L}: {\cal SU}_X(r)\dash |{\cal L}|^*$
associated to
${\cal L}$? Let us observe that we can define a natural (rational) map of
${\cal
SU}_X(r)$ to the linear system $|r\Theta|$, where $\Theta$ denotes the
canonical
Theta divisor on $J^{g-1}(X)$: for $E\in {\cal SU}_X(r)$, define
$$\theta(E):=\{L\in J^{g-1}\ |\ h^0(E\otimes L)\ge 1\}\ .$$
It is easy to see that $\theta(E)$  either is a divisor in $J^{g-1}(X)$ which
belongs to
the linear system  $|r\Theta|$, or is equal to $J^{g-1}(X)$. This last case can
unfortunately occur (see \S 3 below), but only for special $E$'s, so we get a
rational map
$$ \theta: {\cal SU}_X(r)\dash |r\Theta|\ .$$  \medskip

{\pc THEOREM} 2. --  {\it There is a canonical isomorphism $$H^0( {\cal
SU}_X(r),
{\cal L}) \iso H^0(J^{g-1}(X),{\cal O}(r\Theta))^*\ ,$$making the following
diagram
commutative}: \vskip-12pt$$\diagram{
&&|{\cal L}|^*&\cr
{\cal SU}_X(r) &\qquad & \vfl{}{\kern -3pt\wr}{6mm}\cr
&&|r\Theta| &.\cr
\put(15,18){\scriptstyle\varphi^{}_{\scriptscriptstyle\cal L}}
\put(16,4){\scriptstyle\theta}
}$$\vskip-12pt
\ind This  is proved in [B 1] for the rank $2$ case and in  [B-N-R] in
general. Let me say a few words about the proof. For $L$ in $J^{g-1}(X)$ denote
by
$H_L$ the hyperplane in $|r\Theta|$ consisting of divisors passing through $L$.
One
has $\theta^*H_L=\Theta_L$, so we get a linear map $\theta^*:
 H^0(J^{g-1}(X),{\cal O}(r\Theta))^*\longrightarrow H^0( {\cal SU}_X(r),
{\cal L})$ whose transpose  makes the above diagram commutative. It is
easy to show that $\theta^*$ is injective, hence {\it the whole problem is to
prove
that} $\dim  H^0( {\cal SU}_X(r),{\cal L})=r^g$. This was done by constructing
an
$r$\tx to-one covering $\pi:Y\rightarrow X$ such that the pushforward map
$\pi_*:JY\dash {\cal SU}_X(r)$ is dominant, which gives an injective map of
$H^0( {\cal SU}_X(r),{\cal L})$ into $H^0(JY,{\cal O}(r\Theta))$. Note that
surjectivity
of $\theta^*$ means that the linear system $|{\cal L}|$ is spanned by the
divisors
$\Theta_L$ for $L$ in $J^{g-1}(X)$.
 \ind This theorem provides  a relatively concrete
description of the map $\varphi^{}_{\cal L}$, and gives some hope of being able
to
analyze the nature of this map, e.g. whether it is a morphism, an embedding,
and so
on. As we will see, this is a rather intriguing question, which is far from
being
completely understood. We first consider whether this map is everywhere defined
or not.
\vskip1cm {\bf 3. Base points}\smallskip

\ind It follows  from thm. 2 (more precisely, from the fact that the divisors
$\Theta_L$ span the linear system $|{\cal L}|$) that the base points of
$|{\cal L}|$ are
the elements $E$ of $ {\cal SU}_X(r)$ such that $\theta(E)=J^{g-1}(X)$, that is
$H^0(E\otimes L)\not=0$ for {\it all} line bundles $L$ of degree $g-1$. The
existence
of such vector bundles has been first observed by Raynaud [R]. Let me summarize
his results in our language:
\medskip {\pc THEOREM} 3. -- a) {\it For $r=2$, the linear system $|{\cal L}|$
has no
base points.}
\ind b) {\it For $r=3$, $|{\cal L}|$ has no base points if  $g=2$, or if $g\ge
3$ and $X$ is generic.}
\ind c) {\it Let $n$ be an integer $\ge 2$ dividing $g$. For $r=n^g$, the
system
$|{\cal L}|$ has
base points.}
\ind In case c), Raynaud's construction gives only finitely many base points.
This
leaves open a number of questions, which I will regroup under the same heading:
\smallskip
{\bf Q 1}. --  {\it Can one find more examples {\rm (}e.g. for other values of}
$r$)? {\it Is the base locus of dimension} $>0$?  {\it On the opposite side,
can one find a reasonable bound on the dimension of the base locus}?
\smallskip \ind Since the linear system $|{\cal L}|$ has (or may have...) base
points, we have to turn to its multiples. Here we have the following result of
Le Potier [LP], improving an idea of [F 1]:
\smallskip   {\pc PROPOSITION} 1.-- {\it For $\displaystyle k>{1\over
4}r^3(g-1)$, the linear system  $|{\cal L}^k|$ is base-point free.}\smallskip
\ind In fact Le Potier proves a slightly stronger statement: {\it given
$E\in{\cal SU}_X(r)$ and $\displaystyle k>{1\over 4}r^3(g-1)$, there exists
a vector bundle $F$ on $X$ of rank $k$ and degree $k(g-1)$ such that
$H^0(X,E\otimes F)=0$} (in other words,  what may fail with a \lb e always
works  with a rank $k$ vector bundle). Then $\Theta_F:=\{E\in{\cal
SU}_X(r)\ |\ H^0(X,E\otimes F)\not=0\}$ is a divisor of the linear system
$|{\cal L}^k|$ which does not pass through $E$, hence the proposition.
\smallskip \ind   The bound on $k$ is certainly far from optimal; in view of
prop. 1 the most optimistic guess is\smallskip
 {\bf Q 2}. -- {\it Is $|{\cal L}^2|$ base-point free}?
\medskip
\ind Let me also mention the following question of Raynaud [R]:
\smallskip
{\bf Q 3}. -- {\it Given  $E\in{\cal SU}_X(r)$, does there exist an \'etale
covering $\pi:Y\rightarrow X$ such that} $\theta(\pi^*E)\not=J^{g-1}(Y)$?
\vskip1cm {\bf 4. Rank 2}\smallskip
\ind The rank $2$ case is of course the
simplest one; it has two special features. On one hand, by   thm. 3 a) (which
is  quite easy), we know that in this case $\varphi^{}_{\cal L}$ is a {\it
morphism}; we also know
  that this morphism is finite because
${\cal L}$ is ample. On the other hand,  the linear system
$|2\Theta|$ on $J^{g-1}(X)$ is particularly interesting because it contains
the {\it Kummer variety}  ${\cal K}_X$ of $X$. Recall that ${\cal K}_X$ is
the quotient of the Jacobian $JX$ by the involution $a\mapsto -a$, and that
the map\note{1}{As usual $\Theta_a$ denotes the translate of $\Theta$ by
$a$.} $a\mapsto \Theta_a+\Theta_{-a}$ of $JX$ to $|2\Theta|$ factors
through an embedding $\kappa:{\cal K}_X\mono |2\Theta|$. The non-stable
part of ${\cal SU}_X(2)$
 consists of  \vb es of the form $L\oplus L^{-1}$, for $L$ in $JX$, and can
therefore be identified with ${\cal K}_X$; recall that for $g\ge 3$ this is the
singular locus of ${\cal SU}_X(2)$. Thm. 2 thus gives the following
commutative diagram\vskip-8pt$$\diagram{ &&|{\cal L}|^*\cr {\cal
K}\mono{\cal SU}_X(2) &\qquad & \vfl{}{\kern -3pt\wr}{8mm}\cr
&&|2\Theta| &.\cr
\fleche(23,15)\dir(3,2)\long{12}
\fleche(23,13)\dir(3,-2)\long{12}
\fleche(23,15)\dir(3,2)\long{12}
\fleche(4,13)\dir(3,-1)\long{29}
\put(28,22){\scriptstyle\varphi_{\scriptscriptstyle\cal L}}
\put(29,11){\scriptstyle\theta}
\put(12,7){\scriptstyle \kappa}
}$$\vskip-12pt
\ind  Let me summarize what is known  about the structure of
$\varphi^{}_{\cal L}$ (or, what amounts to the same, of $\theta$).
Remember that the dimension of $|2\Theta|$ is $2^g-1$. \medskip {\pc
THEOREM} 4. -- a) {\it For} $g=2$,
 $\theta$ {\it is  an
isomorphism of ${\cal SU}_X(2)$ onto} $|2\Theta|\cong  {\bf P}^3$
[N-R~1].
\ind  b) {\it For  $g\ge 3$,  $C$ hyperelliptic,
  $\theta$ is $2$-to-$1$ onto a  subvariety of $|2\Theta|$ which can be
described in an explicit way} [D-R].
\ind  c) {\it For $g\ge 3$, $C$ not hyperelliptic, $\theta$ is of degree
one onto its image} [B 1]. {\it Moreover if $g(C)=3$ or if $C$ is} generic,
$\theta$ {\it is an embedding} ([N-R 2];[L],[B-V]).
\ind The genus $3$ (non hyperelliptic) case deserves a special mention:  in
this case Narasimhan and Ramanan prove that $\theta$ is an isomorphism of
${\cal SU}_X(2)$ onto a quartic hypersurface ${\cal Q}_4$ in $|2\Theta|$
$(\cong  {\bf P}^7)$. By the above remark, this quartic is singular along the
Kummer variety ${\cal K}_X$. Now it had been observed long time ago
by Coble [C] that there exists a {\it unique}  quartic hypersurface in
$|2\Theta|$  passing doubly through ${\cal K}_X$! Therefore ${\cal Q}_4$
is nothing but Coble's hypersurface.  \ind Part c) of the theorem leaves open
an obvious question:
\smallskip {\bf Q 4}. -- {\it Is $\theta$  always an embedding for
 $C$ non hyperelliptic}?
\smallskip
\ind The case of a generic curve was proved first by Laszlo [L];  Brivio and
Verra
have recently developed  a more geometric approach [B-V], which might hopefully
lead to a complete answer to the question  -- though  some serious
technical difficulties remain at this moment.
\ind Laszlo's method  is to
look at the canonical maps $$\mu_k:S^kH^0({\cal SU}_X(2),{\cal L})
\longrightarrow H^0({\cal SU}_X(2),{\cal L}^k)\ .$$ Since we know that
${\cal L}^k$ is very ample for some (unknown!) integer $k$, surjectivity of
$\mu_k$ for $k$ large  enough would imply that ${\cal L}$ itself is very
ample. We can even dream of getting surjectivity for all $k$, which would
mean that the image of ${\cal SU}_X(2)$ in $|2\Theta|$ is projectively
normal. In [B 2]
 the situation is completely analyzed for $\mu_2$. Recall that a  {\it
vanishing thetanull} on $X$ can be defined as a line bundle $L$ on $X$ with
$L^{\otimes 2}\cong \omega_X$ and $h^0(L)$ even $\ge 2$ -- this means
that the corresponding theta function on $JX$ vanishes at the origin, hence
the name. Such a line bundle exists only on a special curve (more precisely on
a divisor in the moduli space of curves). Then:
 \medskip {\pc PROPOSITION } 2.-- {\it If $X$ has no vanishing thetanull, the
map $\mu_2$ is an isomorphism of $S^2H^0({\cal SU}_X(2),{\cal L})$ onto }
$H^0({\cal SU}_X(2),{\cal L}^2)$.
\ind More generally, if $X$  has $v$ vanishing thetanulls, one has $\dim\Ker
\mu_2=$\break$\dim \Coker \mu_2=v$. This is only half encouraging  (it shows
that  ${\cal SU}_X(2)$ is {\it not} projectively normal for curves with
vanishing thetanulls), but note that the case $k=2$ should be somehow the
most difficult. On the positive side we have the following results:
\smallskip
{\pc PROPOSITION } 3.--  a) {\it If $X$ has no vanishing thetanull, the
map $\mu_4:S^4H^0({\cal SU}_X(2),{\cal L})$ $\longrightarrow
H^0({\cal SU}_X(2),{\cal L}^4)$ is surjective} [vG-P].
\ind b) {\it If $X$ is} generic, {\it the map $\mu_{k}: S^{k}H^0({\cal
SU}_X(2),{\cal L})\longrightarrow H^0({\cal SU}_X(2),{\cal L}^{k})$ is
surjective for $k$ even $\ge 2g-4$} [L].
 \ind As
already mentioned, b) implies that $\varphi^{}_{\cal L}$ is an embedding
for generic $X$.

\vskip1cm {\bf 5. The Verlinde formula}\smallskip
\ind Trying to understand the maps $\mu_k$ raises inevitably the
question of the
dimension of the
spaces $H^0({\cal SU}_X(r),{\cal L}^k)$. We have seen that even the case
$k=1$ is far from trivial -- this is the essential part of [B-N-R].  So it came
as
a great surprise when the mathematical physicists claimed to have a general
(and remarkable) formula for $\dim H^0({\cal SU}_X(r),{\cal L}^k)$, called
the Verlinde formula [V] (there is actually a more general formula for the
moduli space of principal bundles under a semi-simple  group, but we will stick
to the case of  ${\cal SU}_X(r)$):  \medskip
{\pc THEOREM} 5. -- $\displaystyle \dim H^0({\cal SU}_X(r),{\cal
L}^k)=\Bigl({r\over r+k}\Bigr)^g \sum_{{S\amalg T= [1,r+k]}\atop |S|=r}\
\prod_{ {\scriptscriptstyle s\in S}\atop {\scriptscriptstyle t\in T}}\big |
2\sin \pi{s-t\over  r+k}\,\bigr |^{g-1}$.
\medskip
\ind This form of the formula (shown to me by D. Zagier) is the simplest for
an arbitrary rank; for small $r$ or $k$ I leave as a pleasant exercise to the
reader to simplify it (hint: use
$\displaystyle \prod_{p=1}^{n-1}(2\sin{p\pi\over n})=n$) . One gets $r^g$ in
the
case $k=1$, thus confirming thm. 2, and in the rank $2$ case:
\medskip

{\pc COROLLARY}. -- $\displaystyle \dim H^0({\cal SU}_X(2),{\cal
L}^k)=({k\over 2}+1)^{g-1}\sum_{i=1}^{k+1}{1\over  (\sin{i\pi\over
k+2})^{2g-2}}\ \cdotp$ \medskip
\ind Note that  the spaces $ H^i({\cal SU}_X(r),{\cal
L}^k)$ vanish for $i>0$, by the Kodaira vanishing theorem (or rather its
extension by Grauert and Riemenschneider), since the canonical bundle of
${\cal SU}_X(r)$ is equal to ${\cal L}^{-2r}$ [D-N]. Hence thm. 5 gives
actually $\chi({\cal L}^k)$. The right hand side must therefore be a
polynomial in $k$, and take integral values, which is certainly not apparent
from the formula!  In fact I  know no direct proof of these properties, except
in the case $r=2$.
 \ind The leading coefficient of this polynomial is  ${c_1({\cal L})^n\over
n!}$, where $n=(r^2-1)(g-1)$ is the dimension of ${\cal SU}_X(r)$. This
number, which is the volume of ${\cal SU}_X(r)$ for any K\"ahler metric with
K\"ahler class $c_1({\cal L})$,  has been computed in a beautiful way by
Witten [W 1], using the properties of the Reidemeister torsion of a flat
connection. The result is $${c_1({\cal L})^n\over n!}=r\,(2\pi)^{-2n}{\rm
Vol}\,(SU(r))^{2g-2}\,\sum_V{1\over (\dim V)^{2g-2}}\ ,$$ where $V$
runs over all irreducible representations of $SU(r)$, and the volume of
$SU(r)$ is computed with respect to a
 suitably normalized Haar measure. One should be able to deduce this
formula from thm. 5, but I don't know how to do that except in rank $2$.
\vskip1cm
{\bf 6. The Verlinde formula: finite-dimensional proofs}\smallskip
\ind As soon as the Verlinde formula has been known to mathematicians it
has become a challenge for them to give a rigorous proof, so a wealth of
proofs has appeared in the last few years. I will try to describe the ones I am
aware of. The basic distinction is between the proofs using standard algebraic
geometry, which up to now work only in the case $r=2$, and the proofs which
use infinite-dimensional algebraic geometry to mimic the heuristic approach
of the physicists -- these work for all $r$. Let me start with the
``finite-dimensional" proofs. \smallskip \ind The first proof of this kind is
due
to Bertram and Szenes [B-S]; they use the explicit description of the moduli
space in the hyperelliptic case ({\it cf.} thm.  4 b)) to compute $\chi({\cal
L}^k)$ -- which is the same for all smooth curves (actually they
work with the moduli space ${\cal SU}_X(2,1)$ of \vb es of rank $2$ and
fixed determinant of degree $1$, which has the advantage of being smooth;
they show that $\chi({\cal L}^k)$ is the Euler-Poincar\'e characteristic of a
certain vector bundle ${\cal E}_k$ on ${\cal SU}_X(2,1)$).
 \ind A more instructive proof has been obtained by Thaddeus [T], building
on ideas of Bertram and Bradlow-Daskalopoulos. The idea is to look at pairs
consisting of a rank $2$ \vb e $E$ of fixed, sufficiently high degree, say
$2d$, together with a nonzero section $s$ of $E$. There is a notion of
stability for these pairs, in fact  there are various such notions, depending
on
an integer $i$ with $0\le i\le d$. For each of these values one gets a
 moduli space $M_i$, which is  projective and smooth; the key point is that
one passes from $M_{i-1}$ to $M_i$ (for $i\ge 2$)  by a very simple
procedure called a {\it flip} -- blowing up a smooth subvariety and blowing
down the exceptional divisor in another direction. Moreover $M_0$ is just a
projective space, $M_1$ is obtained by blowing up a smooth subvariety in
$M_0$, while $M_{d-1}$ maps surjectively to  ${\cal SU}_X(2)$. In short
one gets the following diagram $$
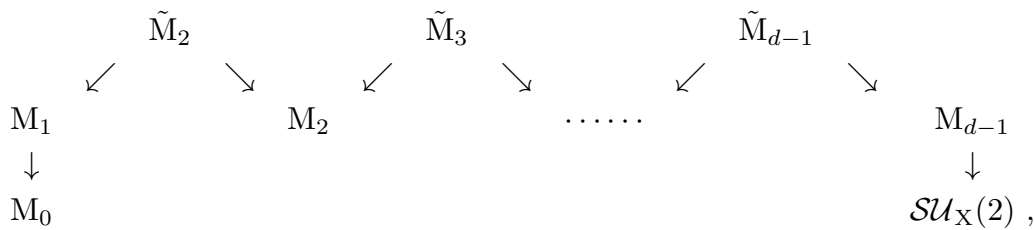
\diagram{
&&\tilde M_2&&&&\tilde M_3&&&&\tilde M_{d-1}\cr
&\swarrow&&\searrow&&\swarrow&&\searrow&&\swarrow&&\searrow\cr
M_1&&&&M_2&&&&\cdots\cdots&&&&M_{d-1}\cr
\downarrow &&&&&&&&&&&&\downarrow\cr
M_0&&&&&&&&&&&&{\cal SU}_X(2)\ ,}$$
from which one deduces (with some highly non-trivial computations) the
Verlinde formula.
\ind Another completely different proof has been obtained
by Zagier (unpublished). Not surprisingly, it is purely computational. Building
on the work of  Atiyah-Bott, Mumford and others, Zagier gives a complete
description of the cohomology ring of ${\cal SU}_X(2,1)$; he is then able to
write down explicitely the Riemann-Roch formula for  $\chi({\cal E}_k)$ (see
above).
\ind Another approach, due, I believe, to Donaldson and
Witten,  starts from Witten's formula for the volume of ${\cal
SU}_X(r)$ (\S 5). More precisely, Witten gives also a formula for the volume
of the moduli space of stable {\it parabolic} bundles; here again the stability
depends  on certain rational numbers, so we get a collection of volumes
indexed by these rational numbers. It turns out that one can recover from
these volumes all the coefficients of the polynomial $\chi({\cal L}^k)$.
\ind The last approach of this type I'd like to mention  has
been developed in [N-Rs] and [D-W 1], and carried out successfully in [D-W
2]. These authors attempt to prove  directly that the spaces  $H^0({\cal
SU}_X(2),{\cal L}^k)$ and the analogous spaces  defined using  parabolic \vb
es obey the so-called {\it factorization rules} (see below).
  Though this approach looks quite promising, the details are
unfortunately quite technical, and an extension to higher rank seems out of
reach.  \vskip1cm {\bf 7. The Verlinde formula: infinite-dimensional
proofs}\smallskip \ind The idea here is to translate in algebro-geometric
terms the methods of the physicists. Actually what the physicists are
interested in is a vector space which plays a central role in Conformal Field
Theory,  the  {\it space of conformal blocks} $B^k_X(r)$. This is defined as
follows: let ${\bf C}((z))$ be the field of formal Laurent series in one
variable.
There is a canonical representation $V_k$ of the Lie algebra ${\goth s\goth
l}_r\bigl({\bf C}((z))\bigr)$ (more precisely, of its universal central
extension), called the {\it basic representation of level $k$}. Let $p\in X$;
the affine algebra $A_X:={\cal O}(X\moins p)$ embeds into ${\bf C}((z))$
(by associating to a function its Laurent expansion at $p$). Then  \smallskip
$$B_X^k(r):= \{\ell\in V_k^*\ |\ \ell(Mv)=0 \quad{\rm for\ all}\quad
M\in{\goth s\goth l}_r(A_X)\ ,\ v\in V_k\}\ .$$  \medskip
 {\pc THEOREM} 6. -- a) {\it There is a canonical isomorphism
$H^0({\cal SU}_X(r),{\cal L}^k)\iso B_X^k(r)$.}
\ind b) {\it The dimension of both spaces is given by the Verlinde formula}
(thm. 5). \ind There are by now several available proofs of these results.
The fact that the dimension of $B_X^k(r)$ is given by the Verlinde formula
follows from the work of Tsuchiya, Ueno and Yamada
 [T-U-Y]. They show that the dimension of $B_X^k(r)$ is independent of the
curve $X$, {\it even if  $X$ is allowed to have double points}. Then it is not
too difficult to express $B_X^k(r)$ in terms of analogous spaces for the
normalization of $X$ (this is called the {\it factorization rules} by the
physicists). One is thus reduced to the genus $0$ case (with marked points),
that is to a problem in the theory of representations of semi-simple Lie
algebras, which is non-trivial in general (actually I know no proof for the
case
of an arbitrary semi-simple Lie algebra), but rather easy for the case of
${\goth s\goth l}_r({\bf C})$.
 \ind Part a)  is proved (independently) in  [B-L] and [F 2]; actually Faltings
proves both a) and b). He considers a smooth curve $X$ degenerating to a
stable curve $X_s$. It is not too difficult to show that $H^0({\cal
SU}_X(r),{\cal L}^k)$ embeds into $B_X^k(r)$, and that the $B_X^k(r)$'s
are semi-continuous so that  $\dim B_X^k(r)\le \dim B_{X_s}^k(r)$.
Therefore the heart of [F 2] is the proof of the inequality  \vskip-12pt$$\dim
B_{X_s}^k(r)\le \dim H^0({\cal SU}_X(r),{\cal L}^k)\ ,$$ on which I cannot
say much, since I don't really understand it (note that  the proof, as well as
that of [T-U-Y], works in the more general set-up of principal bundles).
\smallskip  \ind I would like
to explain in a few words how we construct the  isomorphism\break
$H^0({\cal SU}_X(r),{\cal L}^k)\iso B_X^k(r)$, because I believe its
importance goes far beyond the Verlinde formula. The basic object in the
proof is not ${\cal SU}_X(r)$, but the {\it moduli stack} ${\cal SL}_X(r)$
parametrizing vector bundles $E$ on $X$ together with a trivialization of
$\bigwedge^rE$. Though it appears at first glance as a rather frightening
object, it is both more natural and easier to work with that the moduli space:
basically, working with the moduli stack eliminates all the artificial problems
of
non-representability due to the fact that \vb es have non-trivial
automorphisms. The proof (which is entirely algebraic) has 3 steps:

\par\hskip.5cm 1) We show that {\it the  moduli stack ${\cal SL}_X(r)$
  is isomorphic to the quotient
stack} $SL_r(A_X)\backslash SL_r\bigl({\bf C}((z))\bigr)/SL_r({\bf
C}[[z]])$. The key point here is that a \vb e with trivial determinant is {\it
algebraically trivial} over $X\moins p$ ({\it Hint}: show that  such a bundle
has always a nowhere vanishing section, and use induction on the rank). We
choose a small disk\note{1}{To avoid convergence problems we actually take
$D=\Sp({\cal O})$, where ${\cal O}$ is the completed local ring of $X$ at
$p$, but this makes essentially no difference.} $D\i X$ around $p$, and
consider triples $(E,\rho,\sigma)$ where $E$ is a \vb e on $X$, $\rho$ an
{\it algebraic}  trivialization of $E$  over $X\moins p$ and $\sigma$ a
trivialization of $E$ over $D$.  Over $D\moins p$ these two trivializations
differ by a holomorphic map   $D\moins p \longrightarrow GL_r({\bf C})$
which is meromorphic at $p$, that is given by a Laurent series  $\gamma\in
GL_r\bigl({\bf C}((z))\bigr)$. Conversely given such a matrix $\gamma$ one
can use it to glue together the trivial bundles on $X\moins p$ and $D$ and
recover the triple $(E,\rho,\sigma)$. Since we want $\gamma$ in
$SL_r\bigl({\bf C}((z))\bigr)$ we impose moreover that $\wedge^r\rho$
and $\wedge^r\sigma$ coincide over $D\moins p$. This gives a bijection of
the set  of triples $(E,\rho,\sigma)$ (up to isomorphism)
  onto $SL_r\bigl({\bf C}((z))\bigr)$.
\ind To get rid of the the trivializations, we have to mod out by the
automorphism group  of the trivial bundle over $D$ and $X\moins p$.  We
get the following diagram: $$\nospacedmath\diagram{ \{
E,\,\rho,\,\sigma\}  & \longleftrightarrow  &  SL_r\bigl({\bf C}((z))\bigr)\cr
\vfl{}{}{6mm}& & \vfl{}{} {6mm}\cr \{E,\rho\} &  \longleftrightarrow &
{\cal Q}:=SL_r\bigl({\bf C}((z))\bigr)/SL_r({\bf C}[[z]])\cr \vfl{}{}{6mm} &
& \vfl{}{\pi}{6mm} \cr \{E\} &  \longleftrightarrow & SL_r(A_X)\backslash
SL_r\bigl({\bf C}((z))\bigr)/SL_r({\bf C}[[z]])\ .}$$
Of course I have only constructed a bijection between the set of isomorphism
classes of \vb es on $X$ with trivial determinant and the set of double
classes\break $SL_r(A_X)\backslash SL_r\bigl({\bf C}((z))\bigr)/SL_r({\bf
C}[[z]])$; with some technical work one shows that the construction actually
gives
an isomorphism of stacks.

 \par\hskip.5cm 2) Recall that if $Q=G/H$ is a
homogeneous space, one associates to any character $\chi:H\rightarrow
{\bf C}^*$ a line bundle $L_\chi$ on $Q$: it is the quotient of the
trivial bundle $G\times {\bf C}$ on $G$ by the action of $H$ defined by
$h(g,\lambda)=(gh,\chi(h)\lambda)$. We apply this to the homogeneous
space ${\cal Q}=SL_r\bigl({\bf C}((z))\bigr)/SL_r({\bf C}[[z]])$ (this
is actually an {\it ind-variety}, i.e. the direct limit of an increasing
sequence of projective varieties). By 1) we have a quotient map
$\pi:{\cal Q}\longrightarrow {\cal SL}_X(r)$. The line bundle
$\pi^*{\cal L}$ does {\it not} admit an action of $SL_r\bigl({\bf
C}((z))\bigr)$, but of a group  $\widehat{SL}_r\bigl({\bf C}((z))\bigr)$
which is a central ${\bf C}^*$\tx extension of $SL_r\bigl({\bf
C}((z))\bigr)$. This extension splits over the subgroup $SL_r({\bf
C}[[z]])$, so that ${\cal Q}$ is isomorphic to $\widehat{SL}_r\bigl({\bf
C}((z))\bigr)/({\bf C}^*\times SL_r({\bf C}[[z]]))$. Then $\pi^*{\cal
L}$ is the line bundle $L_\chi$, where $\chi:{\bf C}^*\times SL_r({\bf
C}[[z]])\longrightarrow {\bf C}^*$ is the first projection.
\par\hskip.5cm 3) A theorem of Kumar and Mathieu    provides an
isomorphism $H^0({\cal Q},L^k_\chi)\cong V_k^*$. From this and the
definition of a quotient stack one can identify $H^0({\cal SL}_X(r),{\cal
L}^k)$ with the subspace  of $V_k^*$ invariant under $SL_r(A_X)$. This
turns out to coincide with the  subspace  of $V_k^*$ invariant
under the Lie algebra ${\goth s\goth l}_r(A_X)$, which is by definition
$B_X^k(r)$. Finally a Hartogs type argument gives  $H^0({\cal
SU}_X(r),{\cal L}^k)\cong H^0({\cal SL}_X(r),{\cal L}^k)$.
\bigskip
\ind To conclude let me observe that all the proofs I have mentioned are
rather indirect, in the sense that they involve either degeneration arguments
or sophisticated computations. The simplicity of the formula itself suggests
the following  question:\smallskip
{\bf Q 5}. -- {\it Can one find a direct proof of thm.}$5$?
\ind What I have in mind is for instance a computation of $\chi({\cal L}^k)$
 by simply applying the Riemann-Roch formula; this requires the
knowledge of the Chern numbers of the moduli space.  In [W 2], Witten
proposes some very general conjectures which should give the required
Chern numbers for  ${\cal SU}_X(r)$: the preprint [S]
sketches how the Verlinde formula follows from these conjectures. Jeffrey
and Kirwan have proved some of the Witten's conjectures,  and  I
understand that they are very close to a proof of the Verlinde formula along
these lines.

 \vskip1cm {\bf 8. The strange duality}\smallskip
\ind Let me denote by ${\cal U}^*_X(k)$ the moduli space of semi-stable vector
bundles of rank $k$ and degree $k(g-1)$ on $X$; it is isomorphic (non
canonically)
to ${\cal U}_X(k)$. A special feature of this moduli space is that it carries a
{\it
canonical} theta divisor $\Theta_k$: set-theoretically one has
$$\Theta_k=\{E\in{\cal U}^*_X(k)\ |\ H^0(X,E)\not=0\}\ .$$ Put ${\cal M}:={\cal
O}(\Theta_k)$. Consider the morphism $\tau^{}_{k,r}:{\cal SU}_X(r)\times
{\cal U}^*_X(k) \longrightarrow {\cal U}^*_X(kr)$ defined by
$\tau^{}_{k,r}(E,F)=E\otimes F$. An easy application of the theorem of the
square
shows that $\tau^*_{k,r}({\cal O}(\Theta_{kr}))$ is isomorphic to $pr_1^*{\cal
L}^k\otimes pr_2^*{\cal M}^r$. Now $\tau_{k,r}^*\Theta_{kr}$ is  the divisor of
a
section of this line bundle, well-defined up to a scalar; by the K\"unneth
theorem we
get a linear map $\vartheta_{k,r}: H^0({\cal SU}_X(r),{\cal
L}^k)\longrightarrow
H^0({\cal U}^*_X(k) ,{\cal M}^r)^*$, well-defined up to a scalar. In this
section I want
to discuss the following conjecture: \medskip{\bf Q 6}. -- {\bf Conjecture}:
{\it
The map $\vartheta_{k,r}$ is an isomorphism.} \ind I heard from this statement
3 or
4 years ago, as being well-known to the physicists. The conjecture is discussed
at
length, and extended to \vb es of arbitrary degree, in [D-T].
 \ind Let me discuss a few arguments in favor of the
conjecture. \ind {\it a}) The case $k=1$ is exactly thm. 2.
\ind {\it b}) The two spaces have the same dimension. To prove this one needs
to
compute  the dimension of $H^0({\cal U}^*_X(k) ,{\cal M}^r)$; this is easy
(assuming
the Verlinde formula!) because the map $\tau^{}_{1,k}: {\cal
SU}_X(k)\times J^{g-1}(X)\longrightarrow {\cal U}^*_X(k)$ is an \'etale
(Galois) covering of degree $k^{2g}$, and $\tau^*_{1,k}({\cal M})\cong
pr_1^*{\cal
L}\otimes pr_2^*{\cal O}^{}_J(k\Theta)$.  Therefore we get  $$\hss\dim
H^0(({\cal
U}^*_X(k) ,{\cal M}^r)=\chi({\cal M}^r)={1\over k^{2g}}\ \chi({\cal L}^r)\
\chi({\cal
O}^{}_J(kr\Theta)) ={r^g\over k^g}\ \dim
 H^0(({\cal SU}_X(k) ,{\cal L}^r)\ .\hss$$
Now a quick look at the formula of thm. 5 shows that $k^{-g}\dim
 H^0(({\cal SU}_X(k) ,{\cal L}^r)$ is symmetric in $k$ and $r$,
which proves our assertion.
\ind {\it c}) Therefore it is enough to prove e.g. the surjectivity of the map
$\vartheta_{k,r}$, which  has the following geometric meaning:
\smallskip
${\bf Q\ 6'}$. -- {\it The linear system $|{\cal L}^k|$ in ${\cal SU}_X(r)$ is
spanned by
the divisors  $\Theta_F$, for $F$ in ${\cal U}^*_X(k)$.}

\ind (Recall from \S 3 that $\Theta_F$ is the locus of \vb es $E\in {\cal
SU}_X(r)$
such that $H^0(X,E\otimes F)\not=0$). As an application, taking \vb es $F$ of
the
form $L_1\oplus\ldots\oplus L_k$ with $L_i\in J^{g-1}(X)$, one deduces from
prop. 3 b) that {\it the conjecture holds for $r=2$, $k$ even $\ge 2g-4$}
(in this way we get the result  for a generic curve only, but using the methods
in \S
7 below I can extend it to every curve).

 \vskip1cm{\bf 9. The projective
connection}\smallskip
\ind So far we have considered the moduli space ${\cal SU}_X(r)$ for a fixed
curve
$X$. What can we say of the vector spaces $H^0({\cal SU}_X(r),{\cal L}^k)$ when
the
curve $X$ is allowed to vary? It is again a remarkable discovery of the
mathematical
physicists that these vector spaces are {\it essentially independent of the
curve.}
To explain this in mathematical terms,   consider a family of (smooth)
curves $(X_t)_{t\in T}$, parametrized  by a variety $T$; for $t$ in $T$, let us
denote
by ${\cal L}_t$ the determinant line bundle on ${\cal SU}_{X_t}(r)$. Then:
\medskip
{\pc THEOREM} 7. -- {\it The linear systems  $|{\cal L}_t^k|$ define a {\rm
flat} projective bundle over $T$.}\smallskip
\ind Here again we have by now a number of proofs for this result. The first
mathematical proof  is due to Hitchin [H], following the method used by
Welters in the rank $1$ case; I understand that Beilinson and Kazhdan had a
similar proof (unpublished). A different approach, inspired by the work of the
physicists, appears in [F 1]. Finally, one of the main ingredients in [T-U-Y]
is the
construction of a flat vector bundle over $T$ whose fibre at $t\in T$ is the
space
of conformal blocks $B_{X_t}^k(r)$ (the curves of the family are required to
have a
marked point $p_t\in X_t$, together with a distinguished tangent vector $v_t\in
T_{p_t}(X_t)$). Thanks to thm. 6 this provides still another construction of
our flat
projective bundle. I have no doubt that all these constructions give the same
object but I must confess that I haven't checked it.

\ind Let us take for $T$ the moduli space ${\cal M}_g$ of curves of genus
$g$ (here again, the correct object to consider is the moduli stack, but let me
ignore this). We get a flat projective bundle over ${\cal M}_g$, which
corresponds to
a projective representation  $$\rho_{r,k}:\Gamma_g\longrightarrow
PGL\bigl(H^0({\cal SU}_X(r),{\cal L}^k)\bigr)$$ of the fundamental group
$\Gamma_g=\pi_1({\cal M}_g, \hbox{\eightpoint X})$. This group, called the {\it
modular group} by the physicists and the {\it mapping class group} by the
topologists, is a fundamental object: it carries all the topology of ${\cal
M}_g$. So a
natural question is \smallskip {\bf Q 7}. -- {\it What is the representation
$\rho_{r,k}$}? \ind This is a rather intriguing question.  In the rank one
case, the
analogue of $H^0({\cal SU}_X(r),{\cal L}^k)$ is the space $V_k$
 of $k$\tx th order theta functions;  the group $\Gamma_g$ acts on $V_k$
through
its quotient\note{1}{With a grain of salt when $k$ is odd.}
${\rm Sp}\,(2g,{\bf Z})$, and this action is explicitly described by the
classical
``transformation formula" for theta functions -- which shows in particular that
the
action factors through a finite quotient of ${\rm Sp}\,(2g,{\bf Z})$.
Using thm. 2 we get an analogous description for arbitrary rank in the case
$k=1$. I
expect that the general case is far more complicated and in particular that
$\rho_{r,k}$ doesn't factor through ${\rm Sp}\,(2g,{\bf Z})$, but I know no
concrete example where this happens.
\ind Conformal Field Theory predicts that the connection should be
(projectively) {\it unitary}. This is one of the remaining challenges for
mathematicians:
\smallskip
{\bf Q 8}. -- {\it Find a flat hermitian metric $H_X$ on $H^0({\cal
SU}_X(r),{\cal L}^k)$} (i.e. such that the image of $\rho_{r,k}$ is
contained in $PU(H_X)$).
\ind Here again the rank one case is well-known, and thm. 2 gives an answer in
the

\vskip1cm {\bf 10. Are there generalized theta functions?}\smallskip
\ind I believe that the above results give  some evidence that the spaces
\break$H^0({\cal SU}_X(r),{\cal L}^k)$ are  non-Abelian analogues of the spaces
$H^0(JX,{\cal O}(k\Theta))$.
There is however one aspect of the picture which is missing so far in higher
rank, namely the
analytic description of the sections of ${\cal O}(k\Theta)$ as holomorphic
functions. Clearly the theory of theta functions  cannot be extended in a
straightforward way, if only because ${\cal SU}_X(r)$ is simply connected.
\ind One possible approach is provided by our description of the moduli stack
as a
double quotient $SL_r(A_X)\backslash SL_r\bigl({\bf C}((z))\bigr)/SL_r({\bf
C}[[z]])$ (\S 7). The pull back of the determinant line bundle ${\cal L}$ to
the
group $\widehat{SL}_r\bigl({\bf C}((z))\bigr)$ (a  ${\bf C}^*$\tx extension of
$SL_r\bigl({\bf C}((z))\bigr)$) is trivial, so we should be able to express
sections
of ${\cal L}^k$ as functions on $\widehat{SL}_r\bigl({\bf C}((z))\bigr)$. This
is done
in [B-L] in one particular case: we prove that the pull back of the divisor
$\Theta_{(g-1)p}$ is the divisor of a certain algebraic function on
$\widehat{SL}_r\bigl({\bf C}((z))\bigr)$ known as the $\tau$
function\note{2}{This
idea appears already in [Br], with a slightly different language.}. However it
is not
clear how to express other sections of ${\cal L}$ (and even less of ${\cal
L}^k$) in
a similar way.
\ind There are other possible ways of describing sections of ${\cal L}^k$
by holomor\-phic functions. One of these, explored by
D.~Bennequin, is to pull back ${\cal L}$ to  $SL_r({\bf C})^g$ by the dominant
map
 $SL_r({\bf C})^g\dash {\cal SU}_X(r)$ which maps a
$g$\tx tuple $(M_1,\ldots,M_g)$ to the flat vector bundle $E(\rho^{}_M)\in
{\cal
SU}_X(r)$ associated to the repre\-sen\-ta\-tion
$\rho^{}_M:\pi_1(X)\longrightarrow SL_r({\bf C})$ with
 $\rho(a_i)=I$, $\rho(b_j)=M_j$ (here
  $(a_1,\ldots,a_g,$ $b_1,\ldots,b_g)$ are
the standard generators of $\pi_1(X)$).
\ind Despite these attempts I am afraid we are still far of having a
satisfactory
theory of generalized theta functions. So I will end
up this survey with a loosely formulated question: \smallskip {\bf Q 9}. --
{\it Is
there a sufficiently simple and flexible way of expressing the elements of
$H^0({\cal SU}_X(r),{\cal L}^k)$ as holomorphic functions}?

 \vskip1.8cm \centerline{ REFERENCES} \bigskip \baselineskip13pt
\def\num#1{\item{\hbox to\parindent{\enskip [#1]\hfill}}}
\parindent=1.5cm
\num{B 1} A. {\pc BEAUVILLE}: {\sl 	Fibr\'es de rang $2$ sur
les courbes, fibr\'e d\'eterminant et fonctions th\^eta.} Bull. Soc. math.
France {\bf 116},
431-448 (1988).
 \smallskip
\num{B 2} A. {\pc BEAUVILLE}: {\sl 	Fibr\'es de rang $2$ sur les courbes,
fibr\'e d\'eterminant et fonctions
th\^eta II.}  Bull. Soc. math. France {\bf 119}, 259-291 (1991).
 \smallskip
\num{B-L} A. {\pc BEAUVILLE}, Y. {\pc LASZLO}: {\sl
Conformal blocks and generalized theta functions.} Comm.  Math.
Phys., to appear.
\smallskip
\num{B-N-R} A. {\pc BEAUVILLE}, M.S. {\pc NARASIMHAN}, S. {\pc RAMANAN}: {\sl
Spectral curves and the generalised theta divisor}. J. reine angew. Math. {\bf
398},
169-179 (1989). \smallskip
\num{B-S} A. {\pc BERTRAM}, A. {\pc SZENES}: {\sl Hilbert polynomials of moduli
spaces of rank $2$ vector bundles II.} Topology {\bf 32}, 599-609 (1993).
\smallskip
\num{B-V} S. {\pc BRIVIO}, A. {\pc VERRA}: {\sl The Theta divisor of ${\cal
SU}_C(2)$
is very ample if $C$ is not hyperelliptic and Noether-Lefschetz general}.
Preprint (1993).\smallskip
 \num{Br} J.-L. {\pc BRYLINSKI}: {\sl Loop groups and
non-commutative theta-functions}. Preprint (1990?).\smallskip
 \num C A. {\pc COBLE}: {\sl Algebraic geometry
and theta functions.} AMS Coll. Publi. {\bf 10}, Providence (1929; third
edition,
1969).
 \smallskip
\num{D-W 1} G. {\pc DASKALOPOULOS}, R. {\pc WENTWORTH}: {\sl Local
degenerations of the moduli space of vector bundles and factorization of
rank $2$ theta functions. I} Math. Ann. {\bf 297}, 417-466 (1993).
\smallskip
\num{D-W 2} G. {\pc DASKALOPOULOS}, R. {\pc WENTWORTH}: {\sl
Factorization of rank $2$ theta functions  II: proof of the Verlinde formula.}
Preprint (1994). \smallskip

\num{D-R} U.V. {\pc DESALE}, S. {\pc RAMANAN}: {\sl Classification of \vb
es of rank $2$ on hyperelliptic curves.} Invent. math. {\bf 38}, 161-185
(1976). \smallskip \num{D-T} R. {\pc DONAGI}, L. {\pc TU}: {\sl Theta
functions for $SL(n)$ versus $GL(n)$.} Preprint (1992). \smallskip
 \num{D-N} J.M.
{\pc DREZET}, M.S. {\pc NARASIMHAN}: {\sl Groupe de Picard des vari\'et\'es de
modules de fibr\'es semi-stables sur les courbes alg\'ebriques.} Invent. math.
{\bf
97}, 53-94 (1989).
 \smallskip

\num {F 1} G. {\pc FALTINGS}: {\sl Stable $G$\tx bundles and projective
connections.}
J.
 Algebraic Geometry {\bf 2}, 507-568 (1993). \smallskip

\num {F 2} G. {\pc FALTINGS}: {\sl A proof for the Verlinde formula.} J.
Algebraic Geometry, to appear. \smallskip

\num{vG-P} B. Van {\pc GEEMEN}, E. {\pc PREVIATO}: {\sl Prym varieties and the
Verlinde formula.} Math. Annalen {\bf 294}, 741-754 (1992).
\smallskip
\num{H} N. {\pc HITCHIN}: {\sl Flat connections and geometric quantization.}
Comm.
Math. Phys. {\bf 131}, 347-380 (1990).\smallskip
 \num{L} Y. {\pc LASZLO}: {\sl A propos de l'espace des modules des fibr\'es de
rang
$2$ sur une courbe.} Math. Annalen, to appear. \smallskip \num{N-R 1} M.S. {\pc
NARASIMHAN}, S. {\pc RAMANAN}: {\sl Moduli of vector bundles on a compact
Riemann surface.} Ann. of Math. {\bf 89}, 19-51 (1969). \smallskip
\num{N-R 2} M.S. {\pc NARASIMHAN}, S. {\pc RAMANAN}: {\sl $2\theta$\tx linear
systems  on Abelian varieties.} Vector bundles on algebraic varieties, 415-427;
Oxford University Press (1987).
\smallskip
\num{N-Rs} M.S. {\pc NARASIMHAN}, T.R. {\pc
RAMADAS}: {\sl Factorisation of generalised  Theta functions I.} Invent. math.
{\bf
114}, 565-623 (1993). \smallskip \num{LP} J. {\pc LE} {\pc POTIER}: {\sl
Module des fibr\'es semi-stables et fonctions th\^eta.} Preprint (1993).
\smallskip

\num{R} M. {\pc RAYNAUD}: {\sl Sections des fibr\'es vectoriels sur une
courbe.}
Bull. Soc. math. France {\bf 110}, 103-125 (1982).\smallskip
\num{S} A. {\pc SZENES}: {\sl The combinatorics of the Verlinde formula.}
Preprint (1994).\smallskip \num{T} M.
{\pc THADDEUS}: {\sl Stable pairs, linear systems and the Verlinde formula.}
Preprint (1992). \smallskip
  \num{T-U-Y} A. {\pc TSUCHIYA}, K. {\pc
UENO}, Y. {\pc YAMADA}: {\sl  Conformal field theory on universal family of
stable
curves with gauge symmetries.} Adv. Studies in Pure Math. {\bf 19}, 459-566
(1989).
\smallskip
 \num V E. {\pc VERLINDE}: {\sl Fusion rules and modular transformations
in 2d conformal field theory.} Nuclear Physics {\bf B300}  360-376
(1988). \smallskip
\num {We} A. {\pc WEIL}: {\sl G\'en\'eralisation des fonctions ab\'eliennes.}
J. de Math. P. et App. {\bf 17} (9\up{\`eme} s\'er.), 47-87 (1938).
\smallskip
\num {W 1} E. {\pc WITTEN}: {\sl On quantum gauge theories in two
dimensions.} Comm. Math. Phys. {\bf 141}, 153-209 (1991).
\smallskip
\num {W 2} E. {\pc WITTEN}: {\sl Two-dimensional gauge theories revisited.}
J. Geom. Phys. {\bf 9}, 308-368 (1992).
\vskip1cm
\hfill\hbox to 5cm{\hfill A. Beauville\hfill}

\hfill\hbox to 5cm{\hfill URA 752 du CNRS\hfill}

\hfill\hbox to 5cm{\hfill Math\'ematiques -- B\^at. 425\hfill}

\hfill\hbox to 5cm{\hfill Universit\'e Paris-Sud\hfill}

\hfill\hbox to 5cm{\hfill 91 405 {\pc ORSAY} Cedex, France\hfill}

 \bye